\documentclass[12pt]{article}
\usepackage[dvips]{color}
\usepackage[dvips]{graphicx}
\usepackage{amsfonts}
\definecolor{dgreen}{rgb}{0,.8,.3}
\definecolor{lblue}{rgb}{.2,.3,.7}
\definecolor{red}{rgb}{1.0,0.2,0.2}

\input amssym.def

\input epsf.tex

\begin{document}

\renewcommand{\Box}{\rule{2.2mm}{2.2mm}}
\newcommand{\BOX}{\hfill \Box}

\newtheorem{eg}{Example}[section]
\newtheorem{thm}{Theorem}[section]
\newtheorem{lemma}{Lemma}[section]
\newtheorem{example}{Example}[section]
\newtheorem{remark}{Remark}[section]
\newtheorem{proposition}{Proposition}[section]
\newtheorem{corollary}{Corollary}[section]
\newtheorem{defn}{Definition}[section]
\newtheorem{alg}{Algorithm}[section]
\newtheorem{ass}{Assumption}[section]
\newenvironment{case}
    {\left\{\def\arraystretch{1.2}\hskip-\arraycolsep \array{l@{\quad}l}}
    {\endarray\hskip-\arraycolsep\right.}

\def\argmin{\mathop{\rm argmin}}

\makeatletter
\renewcommand{\theequation}{\thesection.\arabic{equation}}
\@addtoreset{equation}{section} \makeatother

\title{Smoothing
SQP methods
for solving  degenerate nonsmooth  constrained optimization problems with applications to bilevel programs
}

\author{ Mengwei Xu\thanks{\baselineskip 9pt School of Mathematical Sciences, Dalian University of Technology, Dalian 116024, China. E-mail: xumengw@hotmail.com.}, \  Jane J. Ye\thanks{\baselineskip 9pt Department of Mathematics
and Statistics, University of Victoria, Victoria, B.C., Canada V8W 2Y2. E-mail: janeye@uvic.ca.
The research of this author was partially supported by NSERC.} \ and \ Liwei Zhang\thanks{\baselineskip 9pt
School of Mathematical Sciences, Dalian University of Technology, Dalian 116024, China. E-mail: lwzhang@dlut.edu.cn. The research of this author was supported by the National Natural Science
Foundation of China under projects  No. 11071029, No. 91330206 and  No. 91130007.}
}
\date{}
\maketitle

\baselineskip 18pt

{\bf Abstract.} We consider a degenerate nonsmooth and nonconvex optimization problem for which the standard constraint qualification such as the generalized Mangasarian Fromovitz constraint qualification {\rm (GMFCQ)} may not hold.
We use smoothing functions with the gradient consistency property to approximate the nonsmooth functions and  introduce a smoothing sequential quadratic programming (SQP)  algorithm under the exact penalty framework.  We show that any accumulation point of a selected subsequence of the iteration sequence generated by the smoothing SQP algorithm  is a Clarke stationary point, provided that the sequence of multipliers and the sequence of exact penalty parameters are bounded. Furthermore, we propose a new condition called the weakly generalized Mangasarian Fromovitz constraint qualification {\rm  (WGMFCQ)} that is weaker than the GMFCQ. We show that the extended version of the WGMFCQ  guarantees the boundedness of the  sequence of multipliers and the sequence of exact penalty parameters and thus guarantees the global convergence of the smoothing SQP algorithm.  We demonstrate that the WGMFCQ can be satisfied by
bilevel programs for which the GMFCQ never holds.
  Preliminary numerical experiments show that the algorithm is efficient for solving degenerate nonsmooth optimization problem such as  the simple bilevel program.



\newpage

\baselineskip 18pt
\parskip 2pt

\section{Introduction.}
In this paper, we consider the constrained optimization problem of the form
\begin{eqnarray*}
({\rm P})~~~~~~~~~\min && f(x)  \nonumber\\
\rm{s.t.} && g_i(x)\leq 0,\ i=1,\cdots,p,\\
&& h_j(x)=0,\ j=p+1,\cdots,q,
\end{eqnarray*}
where the objective function  and constraint functions
$f, g_i( i=1,\cdots,p), h_j( j=p+1,\cdots,q):\mathbb{R}^n\to \mathbb{R}$ are locally Lipschitz.  In particular, our focus is on solving a degenerate problem for which the generalized Mangasarian Fromovitz constraint qualification {\rm (GMFCQ)} may not hold at a stationary point.

The sequential quadratic programming (SQP) method  is one of the most effective methods for solving smooth constrained optimization problems.
 For the current iteration point $x_k$, the basic idea of the SQP method is to generate a descent direction $d_k$ by solving the following quadratic  programming problem:
\begin{eqnarray*}
\min_d && \nabla f(x_k)^T d+\frac{1}{2} d^T W_k d\\
{\rm s.t.} && g_i(x_k)+\nabla g_i(x_k)^T d\leq 0,\ i=1,\cdots,p,\\
&& h_j(x_k)+\nabla h_j(x_k)^T d= 0,\ j=p+1,\cdots,q,
\end{eqnarray*}
where $\nabla f(x)$ denotes the gradient of function $f$ at $x$ and  $W_k$ is a symmetric positive definite matrix that approximates the Hessian matrix of the Lagrangian function.  Then $d_k$ is used to generate the next iteration point:
$x_{k+1}:=x_k+\alpha_k d_k,$
where the stepsize $\alpha_k$ is chosen to yield a sufficient decrease of a suitable merit function.
The SQP algorithm with $\alpha_k=1$ was first studied by Wilson \cite{wilson}  in which the  exact Hessian matrix of the Lagrangian function is used as $W_k$.   Garcia-Palomares and Mangasarian \cite{gm} proposed to use an estimate to  the Hessian matrix.  Han \cite{han1} proposed  to update the matrix $W_k$ by the  Broyden-Fletcher-Goldfarb-Shanno (BFGS) formula.  
When the stepsize $\alpha_k=1$, the convergence is only local.   To obtain a global  convergence, Han \cite{han2} proposed to use
the classical $l_1$ exact penalty function
as a merit function  to  determine the step size.  While the $l_1$ penalty function is not differentiable,
\cite{py} suggested to use the augmented Lagrange function, which is a smooth function as a merit function.
The inconsistency of the system of the linearized constraints is a serious limitation of the SQP method.  Several techniques  have been introduced to deal with the possible inconsistency.
For example, Pantoja and Mayne \cite{pm} proposed to replace the standard SQP subproblem by the following penalized SQP subproblem:
\begin{eqnarray*}
\min_{d,\xi} && \nabla f(x_k)^T d+\frac{1}{2} d^T W_k d+r_k \xi \\
{\rm s.t.}
&& g_i(x_k)+\nabla g_i(x_k)^T d\leq \xi,\ i=1,\cdots,p,\\
&& -\xi \leq h_j(x_k)+\nabla h_j(x_k)^T d\leq \xi,\ j=p+1,\cdots,q,\\
&& \xi\geq 0,
\end{eqnarray*}
where the penalty parameter $r_k>0$.
Unlike the standard SQP subproblem which may not have feasible solutions,  the penalized SQP subproblem is always feasible for sufficiently large positive constants $r_k$.
%
Other alternative methods for inconsistency of the SQP method are also presented \cite{bh,f,gw,ly,s,tone,zz}.
For nonlinear programs which have some simple bound constraints on some  of the variables, Matthias \cite{m} proposed a projected SQP method which combines the ideas of the projected Newton methods and  the SQP method.


 Recently Curtis and Overton \cite{co} pointed out that applying SQP methods directly to a general  nonsmooth and nonconvex constrained optimization problem will fail in theory and in practices. They employed a process of gradient sampling (GS) method to make the search direction effective in nonsmooth regions and  proved that the iteration points generated by the SQP-GS method converge globally to a stationary point of the exact penalty function  with probability one.
The smoothing method is a well-recognized technique for numerical solution of a nonsmooth optimization problem.  Using a smoothing method, one replaces the nonsmooth function by a suitable smooth approximation,  solves a sequence of smooth problems and  drives the approximation closer and closer to the original problem.   The fundamental question is as follows:  what property a family of the smoothing functions should have  in order for the stationary points of the smoothing problems to approach a stationary point of the original problem? In most of the literature,  a particular smoothing function is employed for the particular problem studied.  It turns out that not all smooth approximations of the nonsmooth function can be used in the smoothing technique to obtain the desired result.
 Zhang and Chen \cite{project}   (see also recent survey on the subject by Chen \cite{chen}) identified  the desired property  as the gradient consistency property.   Zhang and Chen \cite{project} proposed a smoothing projected gradient  algorithm for solving  optimization problems with a convex set constraint by using a family of smoothing functions with  the gradient consistency property to approximate the nonsmooth objective function. They proved that  any accumulation  point of the iteration sequence is
  a Clarke stationary point of  the original nonsmooth  optimization problem.
  Recently \cite{lxy,xy} extended the result of \cite{project} to a class of nonsmooth constrained optimization problem using the projected gradient method and the augmented Lagrangian method respectively. Smoothing functions are proposed  and the  SQP method has been used for the smooth problem  in \cite{fp,jr} to solve
 the mathematical programs with complementarity constraints (MPCC)  and  in  \cite{lqzw,tqw}  to solve  the semi-infinite programming (SIP). In this paper we will combine the SQP method and the smoothing technique to design a smoothing SQP method for a class of general constrained optimization problems with smoothing  functions satisfying  the gradient consistency property. 

For the  SQP method under an 
 exact penalty framework
 to converge globally, usually the set of the  multipliers is required to be bounded (see e.g. \cite{fm}). This  amounts to saying that the MFCQ  is required to hold. For the nonsmooth optimization problem, the corresponding  MFCQ  is referred to as the GMFCQ.
Unfortunately, the GMFCQ is quite strong for certain classes of  problems.  For example, it is well known by now that  the GMFCQ never holds for the bilevel program  \cite{yz}.  Another example of a nonsmooth optimization problem which does not satisfy the GMFCQ is a reformulation of an  SIP \cite{lqzw}. In this paper we propose a new constraint  qualification that is much weaker than the GMFCQ. We call it the weakly  generalized Mangasarian Fromovitz constraint qualification {\rm  (WGMFCQ)}.  WGMFCQ is not a constraint qualification in the classical sense. It is defined in terms of the  smoothing functions and  the  sequence of iteration points  generated by the smoothing algorithm.    In our numerical experiment, WGMFCQ is very easy to satisfy for the bilevel programs.


The rest of the paper is organized as follows.  In Section 2, we present  preliminaries which will be used in this paper and introduce the new constraint  qualification WGMFCQ.  In Section 3, we consider the smoothing approximations of the original problem and propose the smoothing SQP method under an $l_\infty$-exact penalty framework.  Then we establish the global convergence for the algorithm.  In Section 4, we apply the smoothing  SQP method to bilevel programs.  
The final section contains some concluding remarks.

We adopt the following standard notation in this paper. For any two vectors $a$ and $b $ in $\mathbb{R}^n$, we denote  their inner product by $a^T b$. Given a function $G: \mathbb{R}^n\rightarrow \mathbb{R}^m$, we denote its Jacobian by $\nabla G(z)\in \mathbb{R}^{m\times n}$ and, if $m=1$, the gradient $\nabla G(z)\in \mathbb{R}^n$ is considered as a column vector.  
 For a set $\Omega\subseteq \mathbb{R}^n$, we denote   
 the interior, relative interior, the closure, the convex hull, and
  the distance from $x$ to 
 $\Omega$  by
  int $\Omega$, ri $\Omega$, cl $\Omega$, co $\Omega$, and
  dist$(x,\Omega)$ respectively. 
For a matrix $A\in \mathbb{R}^{n\times m}$, $A^T$ denotes its transpose. In addition, we let $\mathbf{N}$ be the set of nonnegative integers and $\exp[z]$ be the exponential function.

\section{Preliminaries and the new constraint qualifications}
In this section, we first present some background materials and results
which will be used later on.  We then discuss the issue of constraint qualification.

Let  $\varphi : \mathbb{R}^n\rightarrow  \mathbb{R}$ be Lipschitz continuous near $\bar{x}$.
The directional derivative of $\varphi$ at $\bar{x}$ in direction $d$ is defined by
\begin{eqnarray*}
\varphi'(\bar{x};d)&:=&\lim_{t\downarrow 0} \frac{\varphi (\bar{x}+td )-\varphi(\bar{x})}{t}.
\end{eqnarray*}
The Clarke generalized directional derivative of $\varphi$ at $\bar{x}$ in direction $d$ is defined by
\begin{eqnarray*}
\varphi^\circ(\bar{x};d):=\limsup_{x\rightarrow\bar{x}, \ t\downarrow 0} \frac{\varphi ({x}+td )-\varphi({x})}{t}.
\end{eqnarray*}
The Clarke generalized gradient of $\varphi$ at $\bar{x}$  is a convex and compact subset of $ \mathbb{R}^n$ defined by
 $$\partial \varphi (\bar{x}):= \{ \xi\in  \mathbb{R}^n: \xi^T d \leq  \varphi^\circ(\bar{x};d), \ \ \forall d \in  \mathbb{R}^n \}.$$
Note that when $\varphi$ is convex, the Clarke generalized gradient coincides with the subdifferential in the sense of convex analysis, i.e.,
 $$\partial \varphi (\bar{x})= \{ \xi\in  \mathbb{R}^n: \xi^T(x-\bar{x}) \leq \varphi (x) -\varphi(\bar{x}), \ \ \forall x \in  \mathbb{R}^n \}$$
and, when $\varphi$ is continuously differentiable at $\bar{x}$, we have $\partial \varphi(\bar{x})=\{\nabla \varphi(\bar{x})\}$.  Detailed discussions of the Clarke generalized gradient and its properties can be found in $\cite{c,clsw}$.



 
 For $\bar x$, a feasible solution of problem $(P)$, we denote by $I(\bar x):=\{i=1,\cdots,p:g_i(\bar x)=0\}$  the active set at $\bar x$.
The following nonsmooth Fritz John type multiplier rule holds  by  Clarke \cite[Theorem 6.1.1]{c}) and the nonsmooth calculus (see e.g. \cite{c}).
\begin{thm}[Fritz John Multiplier Rule]  Let  \(\bar x\) be a local optimal solution of problem $( P)$. Then  there exists $r\geq 0, \lambda_i\geq 0$, $i=1,\cdots,p$, $\lambda_j\in \mathbb{R}$, $j=p+1,\cdots,q$ not all zero such that
\begin{eqnarray}
&& 0\in r\partial f(\bar x)+\sum_{i\in I(\bar{x})} \lambda_i\partial g_i(\bar x)
+\sum_{j=p+1}^q \lambda_j\partial h_j(\bar x).\label{FJ}
\end{eqnarray}
\end{thm}
There are two possible cases in the Fritz John multiplier rule: $r>0$ or $r=0$. Let  $\bar{x}$ be a feasible solution of problem (P). If the Fritz John condition (\ref{FJ}) holds with  $r>0$, then we call  $\bar{x}$ a (Clarke) stationary point of (P).
According to Clarke \cite{c}, any multiplier $\lambda\in \mathbb{R}^q$ with $\lambda_i\geq 0, i=1,\dots, p$ satisfying the Fritz John condition (\ref{FJ}) with $r=0$ is an abnormal multiplier. From the Fritz John multiplier rule, it is easy to see that if there is no nonzero abnormal multiplier then any local optimal solution $\bar{x}$ must be a stationary point. Hence it is natural to define the following constraint qualification.

\begin{defn}[NNAMCQ]
We say that the no nonzero abnormal multiplier constraint qualification $(\rm NNAMCQ)$ holds at a feasible point $\bar x$  of problem $( P)$ if
\begin{eqnarray*}
&& 0\in \sum_{i\in I(\bar{x})} \lambda_i\partial g_i(\bar x)+\sum_{j=p+1}^q \lambda_j\partial h_j(\bar x) \mbox{ and }  \lambda_i\geq 0, \ i\in  I(\bar{x}) \Longrightarrow \lambda_i=0,\lambda_j=0.
\end{eqnarray*}
\end{defn}
It is easy to see that NNAMCQ amounts to saying that any collection of vectors 
$$\{v_i ,i \in I(\bar{x}), v_{p+1},\cdots, v_q\}$$
where  $v_i\in \partial g_i(\bar{x}) (i \in I(\bar{x})), v_j\in \partial h_j(\bar{x}) (j=p+1,\cdots, q)$,
are positively linearly independent.
 NNAMCQ is equivalent to the generalized MFCQ which was first introduced by Hiriart-Urruty \cite{hu}.
\begin{defn}[GMFCQ]
A feasible point $\bar{x}$ is said to satisfy the generalized Mangasarian-Fromovitz constraint qualification ${\rm   (GMFCQ)}$ for problem $( P)$  if  \\
(i) $v_{p+1},\cdots,v_q$ are linearly independent, where  $v_j\in \partial h_j(\bar x)$, $j=p+1,\cdots,q$,\\
(ii) there exists a direction 
$d$ such that
\begin{eqnarray*}
&& v_i^T d<0, \quad \forall v_i\in \partial g_i(\bar x), \ i \in I(\bar{x}),\\
&& v_j^T d=0,\ \quad \forall v_j\in \partial h_j(\bar x) , \ j=p+1,\cdots,q.
\end{eqnarray*}
\end{defn}

In order to accommodate infeasible accumulation points in the numerical algorithm,  we now extend the NNAMCQ and the GMFCQ to allow infeasible points. Note that when $\bar x$ is feasible, ENNAMCQ and EGMFCQ reduce to NNAMCQ and GMFCQ respectively.
\begin{defn}[ENNAMCQ]  We say that the extended no nonzero abnormal multiplier constraint qualification $(\rm ENNAMCQ)$ holds at $ \bar x \in \mathbb{R}^n$
if
\begin{eqnarray*}
&&  0\in  \sum_{i=1}^p \lambda_i\partial g_i(\bar x)+\sum_{j=p+1}^q \lambda_j\partial h_j(\bar x)\mbox{ and }  \lambda_i\geq 0 
,\  i=1,\cdots,p,\\
&& \sum_{i=1}^p \lambda_i g_i(\bar x)+\sum_{j=p+1}^q \lambda_jh_j(\bar x)\geq 0.
\end{eqnarray*}
implies that $ \lambda_i=0,\lambda_j=0$.
\end{defn}
\begin{defn}[EGMFCQ]  A point $\bar{x}\in \mathbb{R}^n$ is said to satisfy the extended generalized Mangasarian Fromovitz constraint qualification {\rm   (EGMFCQ)} for problem {\rm (P)}  if \\
(i) $v_{p+1},\cdots,v_q$ are linearly independent, where $v_j\in \partial h_j(\bar x)$, $j=p+1,\cdots,q$,\\
(ii) there exists a direction 
$d$ such that
\begin{eqnarray*}
&& g_i(\bar x)+v_i^T d<0,\ \quad \forall v_i\in \partial g_i(\bar x), \ i=1,\cdots,p,\\
&& h_j(\bar x)+v_j^T d=0,\ \quad \forall v_j\in \partial h_j(\bar x), \ j=p+1,\cdots,q.
\end{eqnarray*}
\end{defn}
Note that under the extra assumption that the functions $g_i$ are directional differentiable, the EGMFCQ coincides with   the conditions  (B4) and (B5) in \cite{jr}.

Since the set of the Clarke generalized gradient can be large, the ENNAMCQ and the EGMFCQ may be too strong for some problems to hold.   In what follows,  we propose two  conditions that are much weaker than the ENNAMCQ and the EGMFCQ respectively. For this purpose, we first recall the definition of smoothing functions.

\begin{defn}\label{defi3-1} Let $g: \mathbb{R}^n\rightarrow R$ be a locally Lipschitz  function.  Assume that, for a given $\rho>0$, $g_{\rho}: \mathbb{R}^n\rightarrow R$  is   a continuously differentiable function.
We say that $\{g_{\rho}:  \rho>0\}$ is a family of smoothing functions of $g$ if
$\lim\limits_{z\to x,\ \rho\uparrow \infty}g_{\rho}(z)=g(x)$ for any fixed $x\in \mathbb{R}^n$.
\end{defn}

\begin{defn}\label{defi3-2}{\rm \cite{min}} Let $g: \mathbb{R}^n\rightarrow R$ be a locally Lipschitz continuous function. We say that a family of smoothing functions
$\{g_\rho:\rho>0\}$ of $g$ satisfies the gradient consistency property if $\displaystyle \limsup_{z\to x,\, \rho\uparrow \infty}  \nabla g_{\rho}(z)$ is nonempty and
$\displaystyle \limsup_{z\to x,\, \rho\uparrow \infty}  \nabla g_{\rho}(z)\subseteq \partial g(x)$
for any $x\in \mathbb{R}^n$, where $\displaystyle \limsup_{z\to x,\, \rho\uparrow \infty}  \nabla g_{\rho}(z)$ denotes the set of all limiting points
$$\displaystyle \limsup_{z\to x,\, \rho\uparrow \infty}  \nabla g_{\rho}(z):=\Big \{ \lim_{k\rightarrow \infty} \nabla g_{\rho_k}(z_k):
z_k\rightarrow x, \rho_k\uparrow \infty \Big \}.$$
\end{defn}
Note that according to \cite[Theorem 9.61 and Corollary 8.47
(b)]{var}, for a locally Lipschitz function $g$ and its smoothing
family $\{g_\rho: \rho>0\}$, one always has the inclusion $$\partial
g(x) \subseteq \displaystyle co \limsup_{z\to x,\, \rho\uparrow
\infty}  \nabla g_{\rho}(z).$$ Thus our definition of gradient
consistency is equivalent to saying that
$$\partial g(x) = \displaystyle co \limsup_{z\to x,\, \rho\uparrow \infty}  \nabla g_{\rho}(z)$$
which is the definition used in \cite{bhk,chen}.

It is natural to ask if one can always find a family of smoothing functions with the gradient consistency property for a locally Lipschitz function. The answer is yes.
Rockafellar and Wets \cite[Example 7.19 and Theorem 9.67]{var} show that for any locally Lipschitz function $g$, one can construct a family of  smoothing functions  of $g$ with the gradient consistency property by the integral convolution:
\begin{eqnarray*}
g_{\rho}(x):=\int_{\mathbb{R}^n} g(x-y) \phi_{\rho}(y)dy = \int_{\mathbb{R}^n} g(y) \phi_{\rho}(x-y)dy,
\end{eqnarray*}
where $\phi_{\rho}:\mathbb{R}^n\to \mathbb{R}_+$ is a sequence of bounded, measurable functions with $\int_{\mathbb{R}^n} \phi_{\rho}(x)dx=1$ such that the sets $B_{\rho}=\{x:\phi_{\rho}(x)>0\}$ form a bounded sequence converging to $\{0\}$ as $\rho \uparrow \infty$.   Although one can  always generate a family of smoothing functions with the gradient consistency property  by integral-convolution with  bounded supports,  there are many other smoothing functions  which are not generated by the integral-convolution with  bounded supports  \cite{bhk,cc,cm,chen,n}.

Using the smoothing technique,  we approximate the locally Lipschitz
functions $f(x)$, $g_i(x)$, $i=1,\cdots,p$ and $h_j(x)$,
$j=p+1,\cdots,q$ by families of smoothing functions
$\{f_{\rho}(x):\rho>0\}$, $\{g^i_{\rho}(x):\rho>0\}$, $i=1,\cdots,p$
and $\{h^j_{\rho}(x):\rho>0\}$, $j=p+1,\cdots,q$. We also assume
that these families of smoothing functions satisfy the gradient
consistency property. We use certain algorithms to solve the smooth
problem and drive the smoothing parameter $\rho$ to infinity. Based
on the sequence of iteration points of the algorithm, we now define
the new conditions.


\begin{defn}[WNNAMCQ] Let $\{x_k \}$ be a sequence of iteration points for problem
$( P)$  and $\rho_k\uparrow \infty$ as $k\rightarrow \infty$. Suppose that $\bar{x}$ is  a feasible   accumulation point of the sequence $\{x_k \}$.
We say that the weakly no nonzero abnormal multiplier constraint
qualification $(\rm WNNAMCQ)$  based on the smoothing functions
$\{g^i_{\rho}(x):\rho>0\}$, $i=1,\cdots,p$,
$\{h^j_{\rho}(x):\rho>0\}$, $j=p+1,\cdots,q$ holds at
$\bar x$ provided    that
\begin{eqnarray*}
&&  0= \sum_{i\in I(\bar{x})} \lambda_i v_i+\sum_{j=p+1}^q \lambda_j v_j \mbox{ and } \lambda_i\geq 0 ,  \  i\in  I(\bar{x}) \Longrightarrow \lambda_i=0,\lambda_j=0
,\\
\end{eqnarray*} for any $K_0\subset K \subset \mathbf N$ such that $\displaystyle  \lim_{k\rightarrow \infty, k\in K}
x_k=\bar x$ and
\begin{eqnarray*}
 v_i&=& \lim_{k\rightarrow \infty, k\in K_0} \nabla g_{\rho_k}^i (x_k),\  i \in I(\bar{x}),\\
v_j& =& \lim_{k\rightarrow \infty, k\in K_0} \nabla h_{\rho_k}^j
(x_k),\  j=p+1,\cdots,q.
\end{eqnarray*}
\end{defn}
\begin{defn}[WGMFCQ] Let $\{x_k \}$ be a sequence of iteration points for problem
$( P)$  and $\rho_k\uparrow \infty$ as $k\rightarrow \infty$. Let $\bar{x}$ be a feasible
accumulation point  of the sequence $\{x_k \}$.
We say that the weakly generalized Mangasarian Fromovitz constraint
qualification {\rm  (WGMFCQ)} based on the smoothing functions
$\{g^i_{\rho}(x):\rho>0\}$, $i=1,\cdots,p$,
$\{h^j_{\rho}(x):\rho>0\}$, $j=p+1,\cdots,q$ holds at  $\bar{x}$
provided the following conditions hold.  For  any $K_0\subset K \subset \mathbf N$ such that $\displaystyle  \lim_{k\rightarrow \infty, k\in K}
x_k=\bar x$ and any
\begin{eqnarray*}
v_i&=& \lim_{k\rightarrow \infty, k\in K_0} \nabla g_{\rho_k}^i (x_k),\ i\in I(\bar{x})\\
v_j &=& \lim_{k\rightarrow \infty, k\in K_0} \nabla h_{\rho_k}^j
(x_k),\  j=p+1,\cdots,q,
\end{eqnarray*}
(i) $v_{p+1},\cdots,v_q$ are linearly independent;\\
(ii) there exists a direction 
$d$  such that
\begin{eqnarray*}
&& v_i^T d<0,\ {\rm for}\ {\rm all } \ i \in I(\bar{x}),\\
&& v_j^T d=0,\ {\rm for}\ {\rm all } \ j=p+1,\cdots,q.
\end{eqnarray*}
\end{defn}

We now extend the WNNAMCQ and the WGMFCQ to accommodate infeasible points.
\begin{defn}[EWNNAMCQ] Let $\{x_k \}$ be a sequence of iteration points for problem
$( P)$  and $\rho_k\uparrow \infty$ as $k\rightarrow \infty$. Let $\bar{x}$ be a 
accumulation point  of the sequence $\{x_k \}$.
We say that the extended weakly  no nonzero abnormal multiplier
constraint qualification $(\rm EWNNAMCQ)$ based on the  smoothing
functions $\{g^i_{\rho}(x):\rho>0\}$, $i=1,\cdots,p$,
$\{h^j_{\rho}(x):\rho>0\}$, $j=p+1,\cdots,q$ holds at $\bar{x}$ provided that
\begin{eqnarray}
&&  0=  \sum_{i=1}^p \lambda_i v_i+\sum_{j=p+1}^q \lambda_j v_j
\mbox{ and } \lambda_i\geq 0,\  i=1,\cdots,p,\label{nna}\\
&& \sum_{i=1}^p \lambda_i g_i(\bar x)+\sum_{j=p+1}^q
\lambda_jh_j(\bar x)\geq 0.\label{nna2}
\end{eqnarray}
implies that $\lambda_i=0,\lambda_j=0$
 for any $K_0\subset K \subset \mathbf N$ such that $\displaystyle  \lim_{k\rightarrow \infty, k\in K}
x_k=\bar x$  and
\begin{eqnarray*}
v_i&=& \lim_{k\rightarrow \infty, k\in K_0} \nabla g_{\rho_k}^i (x_k),\ i=1,\cdots,p,\\
v_j &=& \lim_{k\rightarrow \infty, k\in K_0} \nabla h_{\rho_k}^j
(x_k),\  j=p+1,\cdots,q.
\end{eqnarray*}

\end{defn}
\begin{defn}[EWGMFCQ] Let $\{x_k \}$ be a sequence of iteration points for problem
$( P)$  and $\rho_k\uparrow \infty$ as $k\rightarrow \infty$.  Let $\bar{x}$ be a 
accumulation point  of the sequence $\{x_k \}$.
We say that the extended weakly  generalized Mangasarian Fromovitz
constraint qualification {\rm  (EWGMFCQ)} based on the smoothing
functions $\{g^i_{\rho}(x):\rho>0\}$, $i=1,\cdots,p$,
$\{h^j_{\rho}(x):\rho>0\}$, $j=p+1,\cdots,q$ holds at $\bar{x}$ provided that the following conditions hold. For any  $K_0\subset K \subset \mathbf N$ such that $\displaystyle  \lim_{k\rightarrow \infty, k\in K}
x_k=\bar x$ and any
\begin{eqnarray*}
v_i&=& \lim_{k\rightarrow \infty, k\in K_0} \nabla g_{\rho_k}^i (x_k),\ i=1,\cdots,p,\\
v_j &=& \lim_{k\rightarrow \infty, k\in K_0} \nabla h_{\rho_k}^i
(x_k),\  j=p+1,\cdots,q,
\end{eqnarray*}
(i) $v_{p+1},\cdots,v_q$ are linearly independent;\\
(ii) there exists a nonzero direction 
$d$ such that
\begin{eqnarray}
&& g_i(\bar x)+v_i^T d<0,\ {\rm for}\ {\rm all } \ i=1,\cdots,p,\label{mf1}\\
&& h_j(\bar x)+v_j^T d=0,\ {\rm for}\ {\rm all } \
j=p+1,\cdots,q.\label{mf2}
\end{eqnarray}
\end{defn}

Due to the gradient consistency property, it is easy to see that
the EWNNAMCQ and the EWGMFCQ are weaker than the ENNAMCQ and the EGMFCQ respectively in general. We finish this section with an equivalence between the EWGMFCQ and EWNNAMCQ.
\begin{thm} \label{thm2.1}
The following implication always holds:
$${\rm EWGMFCQ } \Longleftrightarrow {\rm EWNNAMCQ} .$$
\end{thm}
{\bf Proof.} We first show that EWGMFCQ implies EWNNAMCQ. To the contrary we suppose  that EWGMFCQ holds but  EWNNAMCQ does not hold which means that there exist scalars $\lambda_i\in \mathbb{R}$, $i=1,\cdots,q$ not all zero such that conditions $(\ref{nna})-(\ref{nna2})$ hold.  Suppose that $d$ is the direction that satisfies the condition (ii) of EWGMFCQ.
Due to the  the linear independence of $v_{p+1},\cdots,v_q$ (condition (i) of EWGMFCQ),  the scalars   $\lambda_i, i=1,\dots, p$ can not be all equal to zero.  
Multiplying both
sides of condition $(\ref{nna})$ by $d$, it follows from
conditions $(\ref{mf1})$ and $(\ref{mf2})$ that
\begin{eqnarray*}
 0&=&  \sum_{i=1}^p \lambda_i v_i^T d+\sum_{j=p+1}^q \lambda_j v_j^T d\\
 &<& -\sum_{i=1}^p \lambda_i g_i(\bar x)-\sum_{j=p+1}^q \lambda_j h_j(\bar x)\leq 0,\\
\end{eqnarray*}
which is a contradiction.  Therefore,  EWNNAMCQ holds.

We now prove the reverse implication.  Assume 
the EWNNAMCQ holds.
   EWNNAMCQ implies (i) of EWGMFCQ.  If both (i) and (ii) of EWGMFCQ hold, we are done.
Suppose that the condition (ii) of EWGMFCQ does not hold; that is,
there exists  a subsequence $K_0\subset K\subset N$ and $v_1,\cdots,v_q$ with $\displaystyle  \lim_{k\rightarrow \infty, k\in K}
x_k=\bar x$ and
\begin{eqnarray*}
v_i&=& \lim_{k\rightarrow \infty, k\in K_0} \nabla g_{\rho_k}^i (x_k),\ i=1,\cdots,p,\\
v_j &=& \lim_{k\rightarrow \infty, k\in K_0} \nabla
h_{\rho_k}^j(x_k),\  j=p+1,\cdots,q,
\end{eqnarray*}
 such that $(\ref{mf1})$ and $(\ref{mf2})$ fail to hold.  Let $A:=[v_1,\cdots,v_q]$ be the matrix with $v_1,\dots, v_q$ are columns and
\begin{eqnarray*}
&& S_1:=\{z:z=A^T d, \forall d\},\\
&& S_2:=\{z:z_i<-g_i(\bar x),\ i=1,\cdots,p,\ z_j=-h_j(\bar x),\
j=p+1,\cdots,q\}.
\end{eqnarray*}
Then the convex sets ri $S_1$ and ri ${\rm cl} S_2$ are nonempty and
disjoint.  By the separation theorem, there exists $y\in
\mathbb{R}^q$, $\|y\|\neq 0$ such that $y^T z\geq 0, \forall z\in
S_1$ and   $y^T z\leq 0$, $\forall z\in {\rm cl}S_2$.
By taking $z\in {\rm cl} S_2$ such that  $z_j, j=p+1,\dots, q $ are constants  and $z_i \rightarrow -\infty, i\in\{1,\dots, p\}$, we conclude that
\begin{equation}
y_i\geq 0, \quad i=1, \dots, p.
\label{contra1}
\end{equation}
Choosing $z\in {\rm cl}S_2$ with  $z_i=-g_i(\bar x),\ i=1,\cdots,p,\
z_j=-h_j(\bar x),\ j=p+1,\cdots,q$ we have
\begin{equation}
\sum_{i=1}^p y_i g_i(\bar x)+\sum_{j=p+1}^q y_j h_j(\bar x)=-y^Tz
\geq 0.\label{contra2}
\end{equation}
Select an arbitrary $d$. Then $z=A^T d\in S_1$, $z'=A^T (-d)\in S_1$ and hence
$$\sum_{i=1}^p y_i v_i^Td+\sum_{j=p+1}^q y_j v_j^T d=y^Tz\geq 0,$$
$$\sum_{i=1}^p y_i v_i^T(-d)+\sum_{j=p+1}^q y_j v_j^T (-d)=y^Tz'\geq 0.$$
That is,
\begin{equation}
\sum_{i=1}^p y_i v_i+\sum_{j=p+1}^q y_j v_j=0.\label{contra3}
\end{equation}
From the EWNNAMCQ, conditions (\ref{contra1})-(\ref{contra3})  imply that $y=0$, which is a contradiction.  Thus the condition (ii)
must hold.   The proof is therefore complete.
\BOX

In the case when there is only one inequality constraint and no equality constraints in problem (P), the EWNNAMCQ and EWGMFCQ at   $\bar{x}$ reduces to the following condition:
there is no $K_0\subset K \subset \mathbf N$ such that $\displaystyle  \lim_{k\rightarrow \infty, k\in K}
x_k=\bar x$ and $\displaystyle \lim_{k\rightarrow \infty, k\in K_0} \nabla g_{\rho_k}^1 (x_k)\not =0$. This condition is slightly weaker than a similar condition  \cite[(B4)]{lqzw} which requires that there is no $K_0 \subset \mathbf N$ such that  $\displaystyle \lim_{k\rightarrow \infty, k\in K_0} \nabla g_{\rho_k}^1 (x_k)\not =0$.

\section{Smoothing SQP method}
In this section we design the smoothing SQP algorithm and prove its convergence.

Suppose that  $\{g^i_{\rho}(x):\rho>0\}$ and
$\{h^j_{\rho}(x):\rho>0\}$ are  families of smoothing functions for $g_i, h_j$ respectively. 
Let $x_k$ be the current iterate and $(W_k,r_k,\rho_k)$ be current updates of the positive definite matrix, the penalty parameter and the smoothing parameter respectively. We will try to find a descent direction of a smoothing merit function by using the smoothing SQP subprogram.  In order to overcome the inconsistency of the smoothing SQP subprograms, following Pantoja and Mayne \cite{pm} , we solve the penalized   smoothing SQP subprogram:
\begin{eqnarray*}
({\rm QP})_k~~~~~~\min_{d\in \mathbb{R}^n,\xi\in \mathbb{R}}
&& \nabla f_{\rho_k}(x_k)^T d+\frac{1}{2} d^T W_k d+r_k \xi\\
{\rm s.t.} && g^i_{\rho_k}(x_k)+\nabla g^i_{\rho_k}(x_k)^T d\leq \xi,\ i=1,\cdots,p,\\
&&  h^j_{\rho_k}(x_k)+\nabla h^j_{\rho_k}(x_k)^T d\leq \xi,\ j=p+1,\cdots,q,\\
&&  -h^j_{\rho_k}(x_k)-\nabla h^j_{\rho_k}(x_k)^T d\leq \xi,\ j=p+1,\cdots,q,\\
&& \xi\geq 0.
\end{eqnarray*}
If $(d_{k}, \xi_k)$ is a solution of $({\rm QP})_{k}$, then its Karush-Kuhn-Tucker (KKT) condition can be written as:
\begin{eqnarray}
&&0=\nabla  f_{\rho_k}(x_{k})+W_{k} d_{k}+\sum_{i=1}^p \lambda_{i,k}^g \nabla g^i_{\rho_k}(x_{k})
    + \sum_{j=p+1}^q (\lambda_{j,k}^{+} - \lambda_{j,k}^{-}) \nabla h^j_{\rho_k}(x_{k}),\label{kkt1}\\
&& 0= r_k- \left( \sum_{i=1}^p \lambda_{i,k}^g +
     \sum_{j=p+1}^q(\lambda_{j,k}^{+} + \lambda_{j,k}^{-} ) +\lambda_k^{\xi} \right),    \label{kkt0}\\
 && 0\leq \lambda_{i,k}^g \perp (g^i_{\rho_k}(x_k)+\nabla g^i_{\rho_k}(x_k)^T d_k- \xi_k )\leq 0,
      \ i=1,  \cdots,p,\label{kkt2}\\
 && 0\leq \lambda_{j,k}^{+} \perp (h^j_{\rho_k}(x_k)+\nabla h^j_{\rho_k}(x_k)^T d_k- \xi_k )\leq 0,\
      j=p+1,\cdots,q,\label{kkt3}\\
 && 0\leq \lambda_{j,k}^{-} \perp (-h^j_{\rho_k}(x_k)-\nabla h^j_{\rho_k}(x_k)^T d_k- \xi_k )\leq 0,\
      j=p+1,\cdots,q,\label{kkt4}\\
 && 0\leq \lambda_k^{\xi} \perp -\xi_k \leq 0, \label{kkt5}
\end{eqnarray}
where $\lambda_k=(\lambda^g_k,\lambda^{+}_k,\lambda^{-}_k,\lambda^{\xi}_k)$ is a corresponding Lagrange multiplier.

Let $\rho>0, r>0$. We  define the smoothing merit function by
\begin{eqnarray*}
\theta_{\rho,r}(x):=f_{\rho}(x)
+r \phi_{\rho}(x)
\end{eqnarray*}
where $\phi_{\rho}(x):=\max\{0,g^i_{\rho}(x), i=1,\cdots,p,\  | h_{\rho}^j(x)|,j=p+1,\cdots,q \}$
and  propose the following smoothing SQP algorithm.
\begin{alg}\label{algo3-1}
 Let \( \{\beta,\sigma_{1},\sigma_{2}\}\) be constants in $(0,1)$ with $\sigma_1\leq \sigma_2$, \(\{\sigma,\sigma',\hat{\eta}\}\) be constants in $(1,\infty)$.
 Choose an initial point  \(x_0 \), an initial smoothing parameter $\rho_0>0$, an initial penalty parameter $r_0>0$, an initial   positive definite matrix $W_0\in \mathbb{R}^{n\times n}$ and set \(k:=0\).
\begin{enumerate}
\item  Solve $({\rm QP})_{k}$ to obtain $(d_{k},\xi_{k})$ with the corresponding Lagrange multiplier
$\lambda_{k}=(\lambda^g_{k},\lambda^{+}_{k},\lambda^{-}_{k},\lambda^{\xi}_{k})$, go to Step 2.

\item  If $\xi_{k}=0$, set $r_{k+1}:=r_{k}$ and go to Step 3.  Otherwise, set $r_{k+1}:=\sigma' r_{k}$ and go to Step 3.

\item Let $x_{k+1}:=x_k+\alpha_k d_k$,
where $\alpha_k:=\beta^{l}$, \(l\in \{0,1,2\cdots\}\) is the smallest nonnegative integer satisfying
   \begin{eqnarray}\label{al01}
\theta_{\rho_k,r_{k}}(x_{k+1})-
\theta_{\rho_k,r_{k}}(x_k)
\leq -\sigma_{1}\alpha_k d_{k} W_{k} d_{k}.
   \end{eqnarray}
If
    \begin{eqnarray}\label{al3}
\|d_{k}\|\leq \hat{\eta} \rho_{k}^{-1},
    \end{eqnarray}
   set \(\rho_{k+1}:=\sigma \rho_{k}\) and  go to Step 4.  Otherwise, set \(\rho_{k+1}:= \rho_{k}\) and  go to Step 1.
  In either case, update to  a symmetric positive definite matrix $W_{k+1}$ and $k=k+1$.

\item If a stopping criterion holds, terminate.  Otherwise, go to Step 1.
\end{enumerate}
\end{alg}

We
now show  the global convergence of the smoothing SQP algorithm. For this purpose,
we need  the following standard assumption.
\begin{ass}\label{hess}
There exist two positive constants $m$ and $M$, $m<M$ such that for each $k$ and each $ d\in \mathbb{R}^n$,
$$m\|d\|^2 \leq d^T W_k d \leq M \|d\|^2.$$
\end{ass}

\begin{thm}\label{th3.1}
Suppose that $\{(x_k,\rho_k,d_{k},\xi_k, \lambda_k,r_k, W_k)\}$ is a sequence generated by Algorithm \ref{algo3-1}.
 Then for every $k$,
\begin{equation}
 \theta'_{\rho_k,r_{k}}(x_k, d_{k})\leq -d_{k} W_{k} d_{k} \label{descent}
\end{equation}
and $d_k$ is a descent direction of function $\theta_{\rho_k,r_k}(x)$ at $x_k$ provided by Assumption 3.1 holds.
Furthermore suppose that the Algorithm \ref{algo3-1} does not terminate within finite iterations.  Suppose that
 the sequences $\{x_k\}$ and $\{\lambda_k\}$, $\{r_k\}$ are bounded.
Then $\bar{K}:=\{k:\|d_{k}\|\leq \hat{\eta} \rho_{k}^{-1}\}$ is an infinite set and any accumulation point of sequence $\{x_{k}\}_{\bar{K}}$ is a stationary point of problem $({\rm P})$.
\end{thm}
{\bf Proof.}  
Since $(d_{k}, \xi_k)$ is a solution of $(QP)_{k}$, the KKT conditions $(\ref{kkt1})-(\ref{kkt5})$ hold. 
The directional derivative of the function $x\rightarrow |h^j_{\rho_k}(x)|$ at $x_k$ in direction $d_k$ is
\begin{eqnarray*}
\left\{\begin{array}{cc}
-\nabla h^j_{\rho_k}(x_k)^T d_k,&\ \mbox{if }  \ h^j_{\rho_k}(x_k)<0,\\
|\nabla h^j_{\rho_k}(x_k)^T d_k|,&\ \mbox{if } \ h^j_{\rho_k}(x_k)=0,\\
\nabla h^j_{\rho_k}(x_k)^T d_k,&\ \mbox{if }  \ h^j_{\rho_k}(x_k)>0.
\end{array}\right.
\end{eqnarray*}
Denote the index sets
\begin{eqnarray*}
&&I_k:=\{i=1,\cdots,p: g^i_{\rho_k}(x_k)=\phi_{\rho_k}(x_k)\},\\
&&J_k^+:=\{j=p+1,\cdots,q: h^j_{\rho_k}(x_k)=\phi_{\rho_k}(x_k)\},\\
&&J_k^-:=\{j=p+1,\cdots,q: -h^j_{\rho_k}(x_k)=\phi_{\rho_k}(x_k)\},
\end{eqnarray*}
and $\Gamma_k:=I_k \cup J_k^+\cup J_k^-$.
 Therefore the directional derivative of the function $x\rightarrow \phi_{\rho_k}(x)$ at $x_k$ in direction $d_k$ is
\begin{eqnarray*}
\left\{\begin{array}{cc}
0,\ &\mbox{if }  \ \phi_{\rho_k}(x_k)=0\ \mbox{and}  \ \Gamma_k=\emptyset,\\
\max\{0,\nabla g^i_{\rho_k}(x_k)^T d_k ,i\in I_k,\ |\nabla h^j_{\rho_k}(x_k)^T d_k|, j\in J_k^+\},\ &\mbox{if } \ \phi_{\rho_k}(x_k)=0\ \mbox{and}  \ \Gamma_k\neq \emptyset,\\
\max\{\nabla g^i_{\rho_k}(x_k)^T d_k, i\in I_k,\ \nabla h^j_{\rho_k}(x_k)^T d_k, j\in J_k^+,\\
 -\nabla h^j_{\rho_k}(x_k)^T d_k, j\in J_k^-\},
\ &\mbox{if }  \ \phi_{\rho_k}(x_k)>0.
\end{array}\right.
\end{eqnarray*} 
From $(\ref{kkt2})-(\ref{kkt4})$, we have 
\begin{eqnarray*}
&&\nabla g^i_{\rho_k}(x_k)^T d_k \leq \xi_k - g^i_{\rho_k}(x_k)=\xi_k-\phi_{\rho_k}(x_k),\quad i\in I_k,\\
&&\nabla h^j_{\rho_k}(x_k)^T d_k\leq \xi_k- h^j_{\rho_k}(x_k)=\xi_k-\phi_{\rho_k}(x_k),\quad  j\in J_k^+\\
&&-\nabla h^j_{\rho_k}(x_k)^T d_k\leq \xi_k+ h^j_{\rho_k}(x_k)=\xi_k-\phi_{\rho_k}(x_k),\quad j\in J_k^-.
\end{eqnarray*}
Thus, $\phi_{\rho_k}' (x_k,d_k)\leq \xi_k-\phi_{\rho_k}(x_k)$. Therefore,
\begin{eqnarray*}
 \theta'_{\rho_k,r_{k}}(x_k, d_{k})&=&  \nabla f_{\rho_k}(x_k)^T d_{k}
+r_k  \phi_{\rho_k}' (x_k,d_k)\\
&\leq&  \nabla f_{\rho_k}(x_k)^T d_{k}
+r_k \left(\xi_k-\phi_{\rho_k}(x_k) \right).
\end{eqnarray*}
From $(\ref{kkt0})$  and $(\ref{kkt5})$, we know that if $\xi_k>0$,
\begin{eqnarray*}
 r_k= \left( \sum_{i=1}^p \lambda_{i,k}^g +
     \sum_{j=p+1}^q(\lambda_{j,k}^{+} + \lambda_{j,k}^{-} ) \right),
\end{eqnarray*}
which means
\begin{eqnarray}\label{rlam}
 r_k \xi_k= \left( \sum_{i=1}^p \lambda_{i,k}^g +
     \sum_{j=p+1}^q(\lambda_{j,k}^{+} + \lambda_{j,k}^{-} ) \right)\xi_k.
\end{eqnarray}

By taking conditions $(\ref{kkt1})$, $(\ref{kkt2})-(\ref{kkt4})$ and $(\ref{rlam})$ into account, we obtain that for each $k$,
\begin{eqnarray*}
\lefteqn{ \theta'_{\rho_k,r_k}(x_k, d_k)= \theta'_{\rho_k,r_k}(x_k,d_k)
+ \sum_{i=1}^{p} \lambda_{i,k}^g(g^i_{\rho_k}(x_k)+\nabla g^i_{\rho_k}(x_k)^T d_k-\xi_k)}\\
&&+\sum_{j=p+1}^q  \lambda_{j,k}^{+}(h^j_{\rho_k}(x_k)+\nabla h^j_{\rho_k}(x_k)^T d_k-\xi_k)
+\sum_{j=p+1}^q  \lambda_{j,k}^{-}(-h^j_{\rho_k}(x_k)-\nabla h^j_{\rho_k}(x_k)^T d_k-\xi_k)\\
&\leq & -d_k W_k d_k+\sum_{i=1}^{p} \lambda_{i,k}^g (g^i_{\rho_k}(x_k)-\xi_k)
+\sum_{i=p+1}^{q} \lambda_{j,k}^+ (h^j_{\rho_k}(x_k)-\xi_k)\\
&&+\sum_{i=p+1}^{q} \lambda_{j,k}^- (-h^j_{\rho_k}(x_k)-\xi_k)
+r_k \left(\xi_k-\phi_{\rho_k}(x_k) \right)\\
&\leq& -d_k W_k d_k+r_k \left(\xi_k-\phi_{\rho_k}(x_k) \right)
+\left( \sum_{i=1}^{p} \lambda_{i,k}^g 
+\sum_{i=p+1}^{q} \lambda_{j,k}^+
+\sum_{i=p+1}^{q} \lambda_{j,k}^- \right) (\phi_{\rho_k}(x_k)-\xi_k)\\
&=& -d_k W_k d_k-\left(r_k - \sum_{i=1}^{p} \lambda_{i,k}^g 
-\sum_{i=p+1}^{q} \lambda_{j,k}^+
-\sum_{i=p+1}^{q} \lambda_{j,k}^- \right) \phi_{\rho_k}(x_k)\\
&\leq& -d_k W_k d_k.
\end{eqnarray*}
Hence  the inequality (\ref{descent}) holds.
Since $W_k$ is assumed to be positive definite, it follows that  $d_k$ is a descent direction of function $\theta_{\rho_k,r_k}(x)$ at $x_k$ for every $k$.  Therefore, the algorithm is well-defined.

We now suppose that the Algorithm \ref{algo3-1} does not terminate
within finite iterations. We first prove that there always exists some $d_k$ such that $(\ref{al3})$ holds, thus $\bar{K}$ is an infinite set.  

To the contrary suppose that
$\|d_k\|\geq c_0>0$ for each $k$. Then Assumption 3.1
together with condition $(\ref{al01})$ imply  the existence of a
positive constant $c$ such that $ \theta_{\rho_k,r_k}(x_{k+1})\leq
\theta_{\rho_k,r_k}(x_k)-c$. Consequently,
$(\ref{al3})$ fails. From the boundedness of $\{r_k\}$, we know that
$\xi_k=0$ when $k$ is large.  We can then assume that
there exists  a $\bar k$ large enough such that $\rho_k=\rho_{\bar
k}$ and $r_k=r_{\bar k}$ for $k\geq \bar k$ by the updating rule of
$\rho_k$ and $r_k$.

Since
 the sequence $\{x_k\}$ is bounded, the sequence $\{\theta_{\rho_{\bar k},r_{\bar k}}(x_k)\}$ is bounded below. Moreover $
 \theta_{\rho_k,r_k}(x_{k+1})\leq \theta_{\rho_k,r_k}(x_k)-c
 $, $c>0$, which
  imply that the sequence $\{\theta_{\rho_{\bar k},r_{\bar k}}(x_k)\}$ is monotonously decreasing. Hence we have
\begin{eqnarray*}\label{th21}
\sum_{k\geq \bar{k}} c&\leq&
\sum_{k\geq \bar{k}} \left( \theta_{{\rho_{\bar k}},r_{\bar k}}(x_k)- \theta_{{\rho_{\bar k}},r_{\bar k}}(x_{k+1})\right)\\
&=&  \theta_{{\rho_{\bar k}},r_{\bar k}}(x_{\bar k})- \lim_{k\rightarrow \infty} \theta_{{\rho_{\bar k}},r_{\bar k}}(x_{k})\\
&<& \infty,
\end{eqnarray*}
 which is a contradiction.  Therefore $\bar{K}$ is an infinite set, which also implies that $\rho_k\uparrow\infty$ as $k\to\infty$.

Suppose there exists $K\subseteq \bar{K}$  and $\bar{x}$ such that $\displaystyle\lim_{k\to\infty,k\in K}x_{k}=\bar{x}$.
Since the  sequence $\{\lambda_k\}$  is bounded,
without loss of generality,  assume there exist subsequence $K_1\subset K$ such that $(\lambda^g_k,\lambda^{+}_k,\lambda^{-}_k,\lambda_k^{\xi})\to (\bar{\lambda}^g,\bar{\lambda}^{+},\bar{\lambda}^{-},\bar{\lambda}^{\xi})$ as $k\to\infty, k\in K_1$ and $\bar{\lambda}\geq 0$.
 By the gradient consistency property of $f_{\rho}(\cdot)$, $g^i_\rho(\cdot)$, $i=1,\cdots,p$ and $h^j_\rho(\cdot)$, $j=p+1,\cdots,q$, there exists a subsequence $\tilde{K}_1\subset K_1$ such that
\begin{eqnarray*}
&&\lim_{k\to\infty,\,k\in \tilde{K}_1}\nabla f_{\rho_{{k}}}(x_{k})\in\partial f(\bar{x}),\\
&&\lim_{k\to\infty,\,k\in \tilde{K}_1}\nabla g^i_{\rho_{{k}}}(x_{k})\in\partial g_i(\bar{x}),\ i=1,\cdots,p,\\
&&\lim_{k\to\infty,\,k\in \tilde{K}_1}\nabla h^j_{\rho_{{k}}}(x_{k})\in\partial h_j(\bar{x}),\ j=p+1,\cdots,q.
\end{eqnarray*}
Taking  limits in (\ref{kkt1}) and (\ref{kkt3})-(\ref{kkt5}) as $
k\rightarrow \infty, k\in \tilde{K}_1$, by the gradient consistency
properties and $\xi_{k}\to 0$, it is easy to see that $\bar x$ is a stationary point of
problem $({\rm P})$ and the proof of the theorem is complete. \BOX

%

In the rest of this section, we give a sufficient condition for the boundedness of  sequences $\{r_k\}$ and $\{\lambda_k\}$
. We first give the following result on error bounds.
\begin{lemma}\label{mr}
For each $k\in \mathbf{N}$, $j=1,\cdots,l$, let \(F^j_{k}, F^j: \mathbb{R}^n\to \mathbb{R}\) be   continuously differentiable.  Assume that for each $j=1,\cdots,l$, \(\{F^j_{k}(\cdot)\}\) and \(\{\nabla F^j_{k}(\cdot)\}\) converge to \(F^j(\cdot)\) and \(\nabla F^j(\cdot)\) pointwise respectively  as $k$ goes to infinity.
Let $\hat d$ be the point such that $F^j(\hat d)=0, j=1,\dots, l$.  Suppose that  there exist  $\kappa>0$ and $\delta>0$ such that for all $\mu_j\in [-1,1]$, $j=1,\cdots,l$ not all zero and all $d\in \hat d +\delta B$ it holds that
\begin{eqnarray*}
\left\|\sum_{j=1}^l \mu_j \nabla F^j(d)\right\|>\frac{1}{\kappa}.
\end{eqnarray*}
Then for sufficiently large $k$,
\begin{eqnarray}\label{mr0}
{\rm dist}(\hat d,S_k)\leq \kappa \displaystyle\sum_{j=1}^l |F_k^j(\hat d)|,
\end{eqnarray}
where $
S_k:=\{d\in \mathbb{R}^n:F_k^j(d)=0, j=1,\dots, l\}.$
\end{lemma}
{\bf Proof.}  Denote by  $F(d):=\displaystyle\sum_{j=1}^l |F^j(d)|$, $F_k(d):=\displaystyle\sum_{j=1}^l |F_k^j(d)|$. If $\hat{d}\in S_k$ then (\ref{mr0}) holds trivially.
Now suppose that  $\hat{d} \not \in S_k$.
Since \(F_{k}(\hat{d})\to F(\hat{d})\) as $k\to\infty$, there exists a $\bar k \in \mathbf{N}$ such that $F_k(\hat d)< \kappa^{-1} \delta$ when $k\geq \bar k$.
Let $\varepsilon:=F_k(\hat d)$. Then $\varepsilon \kappa <\delta$.  Take $\lambda\in (\varepsilon \kappa,\delta)$.  Then by Ekeland's  variational principle, there exists an $\omega$ such that $\|\omega-\hat d\|\leq \lambda$, $F_k(\omega)\leq F_k(\hat d)$ and the function $\varphi(d):= F_k(d)+\frac{\varepsilon}{\lambda}\|d-\omega\|$ attains minimum at $\omega$. Hence by the nonsmooth calculus of the Clarke generalized gradient, we have
$$0\in \partial F_k(w)+ \frac{\varepsilon}{\lambda}B$$
where $B$ denotes the closed unit ball of $\mathbb{R}^n$. Thus
$\|v_k\|\leq \frac{\varepsilon}{\lambda}< \frac{1}{\kappa}$,
$\forall v_k\in \partial F_k(\omega)$, for $k\geq \bar k$.  We now show that $F_k(w)=0$ by contradiction. Suppose
that $F_k(w) \not =0$. Then there exists at least one $j$ such that
$F^j_k(w)\not =0$. For such a $j$,  $\partial |F^j_k(w)|=\{\pm
\nabla F^j_k(w)\}$. Therefore there exist $\mu_j^k\in [-1,1]$,
$j=1,\cdots,l$ not all zero  such that $
v_k=\displaystyle\sum_{j=1}^l \mu_j^k \nabla F_k^j(\omega)$.  We
assume that there exist a subsequence $K\subset \mathbf{N}$ and
$\mu_j\in [-1,1]$, $j=1,\cdots,l$ not all zero such that for every
$k\in K$, $F_k(w) \not =0$, $\displaystyle\lim_{k\to \infty,k\in
K}\mu_j^k=\mu_j$, $j=1,\cdots,l$. Since \(\{\nabla F^j_{k}(w)\}_k\)
converge to \(\nabla F^j(w)\), we have $
v:=\displaystyle\lim_{k\to\infty,k\in K}v_k= \sum_{j=1}^l \mu_j
\nabla F^j(\omega)$ and $\|v\|\leq \frac{1}{\kappa}$, which is a
contradiction. The contraction shows that  we must have $F_k(w)=0$ and hence $w\in
S_k$. Therefore we have
\begin{eqnarray*}
{\rm dist}(\hat d,S_k)\leq \|\hat d-\omega\|\leq\lambda.
\end{eqnarray*}
Since this is true for every $\lambda\in (\varepsilon \kappa,\delta)$, we have that for all $k\geq \bar k$
\begin{eqnarray*}
{\rm dist}(\hat d,S_k)\leq \varepsilon \kappa=\kappa |F_k(\hat d)|.
\end{eqnarray*}
\BOX

\begin{thm}\label{th3.3} Assume the Assumption \ref{hess} holds.  Suppose that the Algorithm \ref{algo3-1} does not terminate within finite iterations and  $\{(x_k,\rho_k,d_{k},\xi_k, \lambda_k,r_k)\}$  is a sequence  generated by Algorithm \ref{algo3-1}.   
If the $\rm EWGMFCQ$ holds (or equivalently the $\rm EWNNAMCQ$ holds) at any accumulation point $\bar x$, then the following two statements are true: \\
(a) $\{d_k\}$ and $\{\xi_k\}$ are bounded.\\
(b) $\{r_k\}$ and $\{\lambda_k\}$ are bounded.  Furthermore, when $k$ is large enough, $\xi_k=0$.
\end{thm}
{\bf Proof.}
(a) Assume that there exists a subset $K\subseteq \mathbf{N}$ such that  $\displaystyle\lim_{k\rightarrow \infty, k\in K} x_k=\bar x$.  To the contrary, suppose that $\{d_k\}_K$ is unbounded. Then there exists a subset $K_0\subseteq K$
such that $\displaystyle\lim_{k\rightarrow \infty, k\in {K_0}}\|d_k\| =\infty$ and $\displaystyle \lim_{k\rightarrow \infty,k\in {K_0}} x_k=\bar x$. By the gradient consistency property,
without loss of generality we may assume that
\begin{eqnarray*}
v_i&=& \lim_{k\rightarrow \infty, k\in K_0} \nabla g_{\rho_k}^i (x_k),\ i=1,\cdots,p,\\
v_j &=& \lim_{k\rightarrow \infty, k\in K_0} \nabla h_{\rho_k}^j (x_k),\  j=p+1,\cdots,q,
\end{eqnarray*}
By  EWGMFCQ,   $v_{p+1}, \dots, v_{q}$ are linearly independent and  there exists $\hat d$ such that
\begin{eqnarray*}
&&g_i(\bar x)+v_i^T  \hat d<0,\ i=1,\cdots,p,\\
&&h_j(\bar x)+v_j^T \hat d=0,\ j=p+1,\cdots,q.
\end{eqnarray*}
Since the vectors $\{\displaystyle \lim_{k\to\infty,\,k\in {K_0}}\nabla h^j_{\rho_{{k}}}(x_{k}):  j=p+1,\cdots,q\}$ are linearly independent, it is easy to see that for sufficiently large $k\in {K_0}$, the vectors $\{\nabla h_{\rho_k}^{j}(x_k), j=p+1,\cdots,q\}$ are also linearly independent.  Denote by
\begin{eqnarray*}
 F^j(d)&:=&h_j(\bar x)+v_j^T d,\ j=p+1,\cdots,q,\\
F_k^j(d)&:=&h_{\rho_k}^j(x_k)+\nabla h^j_{\rho_k}(x_k)^T d,\
j=p+1,\cdots,q.
\end{eqnarray*}
 Then  $F^j(\hat d)=0, j=p+1\dots, q$.
 Since $v_{p+1}, \dots, v_{q}$ are linearly independent, there is $\kappa$ such that
 $0<\frac{1}{\kappa}<\min\left\{\displaystyle \left\|\sum_{j=p+1}^q\mu_j v_j\right\|:\mu_j\in [-1,1] \mbox{ not all equal to zero}\right\}$.  By  Lemma \ref{mr},  for sufficient large $k$,
\begin{eqnarray}
{\rm dist}(\hat d,S_k)\leq \kappa \sum_{j=p+1}^q|F_k^j(\hat d)|,\label{errorb}
\end{eqnarray}
where
$
S_k:=\{d\in \mathbb{R}^n: F_k^j(d)=0,\ j=p+1,\cdots,q\}.
$
Since $S_k$ is closed, there exists $\hat{d}_k\in S_k$ such that $\|\hat d-\hat{d}_k\|={\rm dist}(\hat d,S_k)$.
Moreover by virtue of (\ref{errorb}),  the fact that $\displaystyle\lim_{k\rightarrow \infty, k\in {K_0}} F_k^j(\hat d)=F^j(\hat d)=0$ for all $j=p+1,\dots, q$  implies that $\|\hat d-\hat{d}_k\|\to 0$ as $k\to \infty, k\in {K_0}$.   Hence for sufficiently large $k$, we have
\begin{eqnarray}
&& h_{\rho_k}^{j}(x_k)+\nabla h_{\rho_k}^{j}(x_k)^T \hat{d}_k=0, j=p+1,\cdots,q,  \label{eqnh}\\
&& g_{\rho_k}^i(x_k)+\nabla g_{\rho_k}^{i}(x_k)^T \hat{d}_k<0, \ i=1,\cdots,p.\label{eqng}
\end{eqnarray}
(\ref{eqnh})-(\ref{eqng})  imply that $(\hat{d}_k,0)$ is a feasible solution for $(QP)_k$.  Since $(d_k,\xi_k)$ is an optimal solution to problem $(QP)_k$,  we have that for any $k\geq \bar k$, $k\in {K_0}$,
\begin{eqnarray}
\nabla f_{\rho_k}(x_k)^T d_k+\frac{1}{2} d_k^T W_k d_k&\leq &
\nabla f_{\rho_k}(x_k)^T d_k+\frac{1}{2} d_k^T W_k d_k+r_k \xi_k\nonumber\\
&\leq& \nabla f_{\rho_k}(x_k)^T \hat{d}_k+\frac{1}{2} \hat{d}_k^T W_k \hat{d}_k.\label{con1}
\end{eqnarray}
Since $\nabla f_{\rho_k}(x_k)^T \hat{d}_k+\frac{1}{2} \hat{d}_k^T W_k \hat{d}_k $ is bounded, it follows that $\{d_k\}_K$ is bounded from Assumption \ref{hess}.  Since $(d_k,\xi_k)$ are feasible for problem $(QP)_k$, by the definition of the smoothing function and the gradient consistency property, it is easy to see that if $\{d_k\}_K$ is bounded,
then $\{\xi_k\}_K$ is also bounded.  Since $K$ and $\bar x$ are arbitrary subset and arbitrary accumulation point, $\{d_k\}$ and $\{\xi_k\}$ are bounded for the whole sequence.

%

(b) To the contrary, suppose that $\{\lambda_k\}$ is unbounded.
Then there exists a subset ${K_1}\subseteq K$
such that $\displaystyle\lim_{k\rightarrow \infty, k\in {K_1}}\|\lambda_k\| =\infty$ and $\xi_k>0$ for $k\in K_1$ sufficiently large.
 By the gradient consistency property, without loss of generality we may assume that
\begin{eqnarray*}
v_i&=& \lim_{k\rightarrow \infty, k\in K_1} \nabla g_{\rho_k}^i (x_k),\ i=1,\cdots,p,\label{vg}\\
v_j &=& \lim_{k\rightarrow \infty, k\in K_1} \nabla h_{\rho_k}^j (x_k),\  j=p+1,\cdots,q\label{vh},
\end{eqnarray*}
and
$\displaystyle
\lim_{k\to\infty,k\in {{K_1}}}\frac{\lambda_{k} }{\|{\lambda_k}\|}=\bar{\lambda}
$
for some nonzero vector $\bar{\lambda}=(\bar{\lambda}^g,\bar{\lambda}^{+},\bar{\lambda}^{-},\bar{\lambda}^
{\xi})\geq 0$.
Dividing by $\|{\lambda_k}\|$ in both sides of $(\ref{kkt1})$ and letting $k\to \infty$, $k\in {{K_1}}$, we have
\begin{eqnarray}\label{th3.0}
0=\sum_{i=1}^p \bar{\lambda}^g_i  v_i+\sum_{j=p+1}^q (\bar{\lambda}^{+}_j-\bar{\lambda}^{-}_j)  v_j .
\end{eqnarray}

Letting $k\to \infty$, $k\in {{K_1}}$ in conditions $(\ref{kkt2})-(\ref{kkt5})$ and assuming that $(\bar d,\bar{\xi})$ is the limiting point of $\{(d_k,\xi_k)\}_{K_1}$, we have
\begin{eqnarray*}
 && 0\leq \bar{\lambda}_{i}^g \perp (g_i(\bar x)+v_i^T \bar{d}- \bar{\xi} )\leq 0,
      \ i=1,  \cdots,p,\\
 && 0\leq \bar{\lambda}_{j}^{+} \perp (h_j(\bar x)+v_j^T \bar{d}- \bar{\xi} )\leq 0,\
      j=p+1,\cdots,q,\\
 && 0\leq \bar{\lambda}_{j}^{-} \perp (-h_j(\bar x)-v_j^T \bar{d}- \bar{\xi} )\leq 0,\
      j=p+1,\cdots,q,\\
 && 0\leq \bar{\lambda}^{\bar{\xi}} \perp -\bar{\xi}\leq 0.
\end{eqnarray*}
Multiplying both sides of $(\ref{th3.0})$ by $\bar d$, since
\begin{eqnarray*}
 &&  \bar{\lambda}_{i}^g (g_i(\bar x)+v_i^T \bar{d}- \bar{\xi} )= 0,
      \ i=1,  \cdots,p,\\
 && \bar{\lambda}_{j}^{+}  (h_j(\bar x)+v_j^T \bar{d}- \bar{\xi} )= 0,\
      j=p+1,\cdots,q,\\
 &&  \bar{\lambda}_{j}^{-}  (-h_j(\bar x)-v_j^T \bar{d}- \bar{\xi} )= 0,\
      j=p+1,\cdots,q,
\end{eqnarray*}
we have
\begin{eqnarray*}
0&=&\sum_{i=1}^p \bar{\lambda}^g_i  v_i^T \bar{d}+\sum_{j=p+1}^q (\bar{\lambda}^{+}_j-\bar{\lambda}^{-}_j)  v_j^T \bar{d}\\
&= & \sum_{i=1}^p \bar{\lambda}^g_i (\bar{\xi}-g_i(\bar x) ) 
+ \sum_{j=p+1}^q \bar{\lambda}^{+}_j  (\bar{\xi}-h_j(\bar x) )
+ \sum_{j=p+1}^q \bar{\lambda}^{-}_j  (\bar{\xi}+h_j(\bar x) ).
\end{eqnarray*}
Thus,
\begin{eqnarray}\label{th3.1new}
 \sum_{i=1}^p \bar{\lambda}^g_i g_i(\bar x) 
+ \sum_{j=p+1}^q (\bar{\lambda}^{+}_j - \bar{\lambda}^{-}_j)  h_j(\bar x) 
=\sum_{i=1}^p \bar{\lambda}^g_i \bar{\xi}
+ \sum_{j=p+1}^q (\bar{\lambda}^{+}_j + \bar{\lambda}^{-}_j)  \bar{\xi} \geq 0.
\end{eqnarray}
From the EWGMFCQ (equivalently EWNNAMCQ), 
condition (\ref{th3.1new}) together with condition (\ref{th3.0}) imply that $\bar{\lambda}_i^g=0$, $i=1,\cdots,p$ and $\bar{\lambda}^{+}_j - \bar{\lambda}^{-}_j=0$, $j=p+1,\cdots,q$.

Consider the case where $\bar{\lambda}_i^g=0$, $i=1,\cdots,p$ and there exists an index $j\in\{p+1,\cdots,q\}$ such that $\bar{\lambda}^{+}_j = \bar{\lambda}^{-}_j>0$.  Then for sufficiently large $k\in K_1$, $\lambda^{+} _{j,k}>0$ and $\lambda^{-} _{j,k}>0$.  From the complementary condition $(\ref{kkt3})-(\ref{kkt4})$, we must have $\xi_k=0$ for sufficiently large $k\in K_1$, which is a contradiction.

Otherwise, consider the case where $\bar{\lambda}_i^g=0$, $i=1,\cdots,p$ and $\bar{\lambda}^{+}_j = \bar{\lambda}^{-}_j=0$, $j=p+1,\cdots,q$. Then since $\bar{\lambda}$ is a nonzero vector, we must have $\bar{\lambda}^{{\xi}} >0$, which implies that $\lambda^{\xi} _k>0$ for sufficiently large $k\in K_1$.
From the complementarity condition $(\ref{kkt5})$, $\xi_k=0$ for sufficiently large $k\in K_1$, which is a contradiction.

 The contradiction shows that  $\{\lambda_k\}$must be  bounded. By the relationship between $\{\lambda_k\}$ and  $\{r_k\}$, the boundedness of  $\{\lambda_k\}$ implies that boundedness of $\{r_k\}$. 
 Furthermore, from the updating rule of the algorithm, the boundedness of the sequences $\{\lambda_k\}$ and  $\{r_k\}$ implies that when $k$ is large enough, $\xi_k=0$.
We complete the proof.
\BOX

The following corollary follows immediately from Theorems \ref{th3.1} and \ref{th3.3}. 
\begin{corollary}\label{Coro3.1} Let Assumption \ref{hess} hold and suppose that the Algorithm \ref{algo3-1} does not terminate within finite iterations.  Suppose that the sequence $\{x_k\}$ is bounded.
Assume 
the $\rm EWGMFCQ$ (or equivalently {\rm EWNNAMCQ}) holds at any accumulation point of sequence $\{x_k\}$,
then  $\bar{K}:=\{k:\|d_{k}\|\leq \hat{\eta} \rho_{k}^{-1}\}$ is an infinite set and any accumulation point of sequence $\{x_{k}\}_{\bar{K}}$ is a stationary point of problem $({\rm P})$.
\end{corollary}
In the case where the objective function is smooth, there is only one inequality constraint and no equality constraints in problem (P), Corollary  \ref{Coro3.1} extends  \cite[Theorem 4.3]{lqzw}  to allow the general smoothing function instead of the specific smoothing function.

\section{Applications to the bilevel programs}
The purpose of this section is to apply the smoothing  SQP algorithm to 
the bilevel program. 
We illustrate how we can apply our algorithm to solve the bilevel program and we demonstrate through some numerical examples   that although the GMFCQ never holds for bilevel programs, the WGMFCQ may be satisfied easily.

In our numerical experiments, we use the following  method  proposed by Powell \cite{powell2} which is a modification to the  BFGS method for unconstrained optimization problems to update the matrix $W_k$.
  Define $s_k:=x_{k+1}-x_k$ and
  \begin{eqnarray*}
\lefteqn{ y_k:=\nabla f_{\rho_k}(x_{k+1})-\nabla f_{\rho_k}(x_k) - \sum_{i=1}^p \lambda_{i,k}^g (\nabla g_{\rho_k}^i(x_{k+1}) - \nabla g_{\rho_k}^i(x_{k}) )} \\
&&-\sum_{j=p+1}^q (\lambda_{j,k}^+-\lambda_{j,k}^- ) (\nabla h_{\rho_k}^j(x_{k+1}) - \nabla h_{\rho_k}^j(x_{k}) ).
\end{eqnarray*}
The modified $\bar{y}_k$ takes the form
\begin{eqnarray*}
 \bar{y}_k=\left\{\begin{array}{cc}
 y_k,& {\rm if}\ s_k^T y_k\geq 0.2 s_k^T W_k s_k,\\
 \theta_k y_k+(1-\theta_k) W_k s_k,&\ {\rm otherwise},
\end{array} \right.
\end{eqnarray*}
where
$\displaystyle
\theta_k=\frac{0.8 s_k^T W_k s_k}{s_k^T W_k s_k - s_k^T y_k}.
$
We update $W_{k+1}$ by
\begin{eqnarray*}
W_{k+1}=W_k - \frac{W_k s_k s_k^T W_k}{s_k^T W_k s_k} +\frac{\bar{y}_k \bar{y}_k^T}{s_k^T \bar{y}_k}.
\end{eqnarray*}
When the norm of $W_{k+1}$ is too large or too small, e.g.  greater than $10^5$ or smaller than $10^{-5}$, we set $W_{k+1}=I$, where $I$ is the  identity matrix. This way we make sure Assumption 3.1 holds.

In numerical practise, it is impossible to obtain an exact `0', thus we select some small enough $\varepsilon>0$, $\varepsilon'>0$ and change the update rule of $r_k$ and $\rho_k$ to the case when $\xi_k<\varepsilon'$ and 
\begin{eqnarray*}
\|d_{k}\|\leq \max\{\hat{\eta} \rho_{k}^{-1}, \varepsilon\}
    \end{eqnarray*}
respectively.
 We suggest the stopping criterion as follows: for a given $\epsilon_1>0$, we terminate the algorithm at the $k$th iteration if
 $$\|(x^{k-1},y^{k-1})-(x^k,y^k)\|< \epsilon_1.$$

To verify the EWGMFCQ, we consider the following cases. 
When the sequence which generated by the algorithm has more than one accumulation points, we should verify all of the accumulation points.  When the sequence has only one accumulation point (which happens frequently), if the accumulation point is feasible, we verify the WGMFCQ at the point, otherwise we change to another initial point.

In the rest of this section we consider the simple bilevel program
\begin{eqnarray*}
({\rm SBP})~~~~~~\min&& F(x,y)\nonumber\\
{\rm s.t.} && y\in S(x), 
\end{eqnarray*}
where $S(x)$ denotes the set of solutions of the lower level program
\begin{eqnarray*}
({\rm P}_x)~~~~~~~~
\min_{y\in Y} \  f(x,y),
\end{eqnarray*}
where
 $F, f:\mathbb{R}^n\times \mathbb{R}^m \rightarrow \mathbb{R}$ are continuously differentiable and  twice continuously differentiable respectively, and $Y$
 is a compact subset of $\mathbb{R}^m$.  Our smoothing SQP algorithm can easily handle any extra upper level constraint but we omit it for simplicity.  For a general bilevel program, the lower level constraint may depend on the upper level variables. By ``simple'', we mean that the lower level constraint $Y$ is independent of $x$.  Although (SBP) is a simple case of the general bilevel program, it has many applications such as the principal-agent problem \cite{Mirrlees99} in Economics. We refer the reader to \cite{b,d1,d2,sib,vsj} for  applications of  general bilevel programs.

When the lower level program is a convex program in variable $y$, the first order approach to solving a bilevel program is to replace the lower level program by  its KKT conditions.
In the case where $f$ is not convex in variable $y$, Mirrlees \cite{Mirrlees99} showed that this approach may not be  valid in the sense that the true optimal solution for the bilevel problem may not even be a stationary point of the reformulated problem by the first order approach.

For a numerical purpose, Outrata \cite{Outrata} proposed to reformulate a bilevel program as a nonsmooth single level optimization problem by replacing the lower level program by its value function constraint, which in our simple case is
\begin{eqnarray}
({\rm VP})~~~~~~\min && F(x,y)\nonumber\\
{\rm s.t.} && f(x,y)-V(x) = 0,\label{SP}\\
&&x\in  \mathbb{R}^n, y\in  Y\nonumber
\end{eqnarray}
where $\displaystyle V(x) :=\inf_{y\in Y} \ f(x,y)$ is the value function of the lower level problem.
%
Ye and Zhu \cite{yz} pointed out that  the usual constraint qualifications such as the GMFCQ never hold  for problem (VP).
Ye and Zhu \cite{yz,yz1} derived the first order necessary optimality condition for the general bilevel program under the so-called ``partial calmness condition'' under which the difficult constraint (\ref{SP}) is moved to the objective function with a penalty.
 Based on the value function approach, Xu and Ye
\cite{lxy} recently proposed to approximate the value function by its integral entropy function:
\begin{eqnarray*}
\gamma_{\rho}(x)
&:=&-\rho^{-1} \ln\left(\int_{Y} \exp[-\rho f(x,y)]dy\right)
\end{eqnarray*}
and developed a  smoothing projected gradient  algorithm to solve the problem $(\rm VP)$ when the problem (SBP) is partially calm and to solve an approximate bilevel problem $(\rm VP)_\varepsilon$ where the constraint (\ref{SP}) is replaced by $f(x,y)-V(x) \leq  \varepsilon$  for  small $\varepsilon>0$  otherwise.

Unfortunately, the partial calmness condition is rather strong and hence a local optimal solution of a bilevel program may not be a stationary point of  (VP).
Ye and Zhu  \cite{yz2} proposed to study the following combined program by adding  the first order condition of the lower level problem into the problem (VP).
Although the partial calmness condition is a very strong condition for $(\rm VP)$, it is likely to hold for the combined problem under some reasonable conditions \cite{yz2}.

Recently  Xu and Ye \cite{xy} proposed a smoothing augmented Lagrangian method to solve the combined problem with the assumption that each lower level solution lies in the interior of $Y$:
\begin{eqnarray}
({\rm CP})~~~~~~\min_{(x,y) \in \mathbb{R}^n\times Y}  && F(x,y) \nonumber\\
{\rm s.t.} && f(x,y)-V(x) \leq 0, \label{cp1}\\
&& \nabla_y f(x,y)=0,\label{cp2}
\end{eqnarray}
 They showed that if the sequence of penalty parameters is bounded, then any accumulation point is a  Clarke stationary point of (CP).  They argued that since the problem (CP) is very likely to satisfy the partial calmness or the weak calmness condition (see \cite{yz2}), the sequence of penalty parameters is likely to be bounded.


To simplify our discussion so that we can concentrate on the main idea, we make the following assumption
\begin{ass}\label{4-1}
Every optimal solution  of the lower level problem  is an interior point of set $ Y$.
\end{ass} 
In practice, it may be possible to set the set $Y$ large enough so that all optimal solutions of the lower level problem are contained in the interior of $Y$. 
If it is difficult to do so and the set $Y$ can be represented by some equality or inequality constraints then one can use the KKT condition to replace the constraint (\ref{cp2}) in the problem (CP). 

 Since problem (CP) is a nonconvex and nonsmooth optimization problem, in general the best we can do is to  look for its Clarke stationary points.  Since we assume that all lower level solutions lie in the interior of set $Y$, any local optimal solution of (CP) must be the Clarke stationary point of (CP) with the constraint $y\in Y$ removed. Hence the smoothing SQP method introduced in this paper can be used to find the stationary points of (CP).

Let $(\bar{x},\bar{y})$ be a local optimal solution of (CP). Then by  the Fritz John type multiplier rule,  there exist $r\geq 0, \lambda_1\geq 0, \lambda_2\in \mathbb{R}^m$ not all zero such that 
\begin{eqnarray}
&&0\in r\nabla F(\bar{x},\bar{y})+\lambda_1 (\nabla f(\bar x,\bar y)-\partial V(\bar x)\times\{0\} )+\nabla (\nabla_y f)(\bar x,\bar y)^T \lambda_2 
\label{newFJ}
\end{eqnarray}
In the case when $r$ is positive,  $(\bar{x},\bar{y})$ is a stationary point of (CP). A sufficient condition for $r$ to be positive is that $r=0$ in the Fritz John condition in which case $\lambda_1,\lambda_2$ should not be all equal to zero. Unfortunately we now show that $r$ can be always taken as zero in the above Fritz John condition for problem (CP).
Indeed, 
from the definition of $V(x)$,  we always have $f(x,y)-V(x) \geq 0$ for any $y\in Y$. Hence
 any feasible point $(\bar x,\bar y)$ of problem (CP) is always an optimal solution of the problem
$$\displaystyle\min_{(x,y)\in  \mathbb{R}^n\times Y}f( x,y)-V( x), \mbox{ s.t. }\nabla_y f(x,y)=0.$$ By the Fritz John type multiplier rule, there exists $\lambda_1\geq 0, \lambda_2 \in \mathbb{R}^m$ not all equal to zero such that
\begin{eqnarray}
0\in \lambda_1 (\nabla f(\bar x,\bar y)-\partial V(\bar x)\times\{0\} )+\nabla (\nabla_y f)(\bar x,\bar y)^T \lambda_2. 
\label{sp1}
\end{eqnarray}
Observe that (\ref{sp1}) is (\ref{newFJ}) with $r=0$. 
Since $(\lambda_1,\lambda_2)$ is nonzero, we have shown that  the Fritz John condition (\ref{newFJ}) for problem (CP) holds with $r=0$.  In another word,  the NNAMCQ (or equivalently GMFCQ)  for problem (CP) never hold.


 It follows from \cite[Theorem 5.1 and 5.5]{lxy} that the integral entropy function 
 \(\gamma_{\rho}(x) \)  is a smoothing function  with the gradient consistency property for the value function $V(x)$. That is,
$$\displaystyle \lim_{z\to x,\, \rho\uparrow \infty}   \gamma_{\rho}(z)=V(x)\quad \mbox{ and } \quad \emptyset \not =\displaystyle \limsup_{z\to x,\, \rho\uparrow \infty}  \nabla \gamma_{\rho}(z)\subseteq \partial V(x).$$
For a sequence of iteration points $\{(x^k,y^k)\}$, the set $\displaystyle \limsup_{k\rightarrow \infty}\nabla \gamma_{\rho_k}(x^{k}) $ may strictly contain in
 $\partial V(x)$. Therefore while  (\ref{sp1}) holds for some $\lambda_1\geq 0,  \lambda_2 \in \mathbb{R}^m$ not all equal to zero, the following inclusion may  hold only when $\lambda_1=0, \lambda_2=0$:
\begin{eqnarray*}
0\in \lambda_1 ( \nabla f(\bar x,\bar y)-\displaystyle \limsup_{k\rightarrow \infty}\nabla \gamma_{\rho_k}(x^{k})\times\{0\}) +\nabla (\nabla_y f)(\bar x,\bar y)^T \lambda_2 
,
\end{eqnarray*}
And consequently,  the WNNAMCQ may hold. We illustrate this point by using some numerical examples.
In these examples, since $y\in \mathbb{R}$, the problem (CP) has one inequality constraint $f(x,y)-V(x)\leq 0$ and one equality constraint $\nabla_y f(x,y)=0$. Hence the WNNAMCQ
$$0\in \lambda_1 (\nabla f(\bar{x},\bar{y})-\displaystyle \limsup_{k\rightarrow \infty}\nabla \gamma_{\rho_k}(x^{k})\times\{0\}) +\lambda_2\nabla (\nabla_y f)(\bar x,\bar y) ,\quad \lambda_1\geq 0 \Longrightarrow \lambda_1=\lambda_2=0$$
amounts to saying that for any $K_0\subset K\subset {\bf N}$ such that $\displaystyle \lim_{k\rightarrow \infty, k\in K} (x^k,y^k)=(\bar{x},\bar{y})$ and $\displaystyle v=\lim_{k\rightarrow \infty, k\in K_0}\nabla \gamma_{\rho_k}(x^{k})$, the vectors
$$\nabla f(\bar{x},\bar{y})-(v,0) \mbox{ and } \nabla (\nabla_y f)(\bar x,\bar y) $$ are linearly independent.

\begin{eg}\label{Mirr}\rm \cite{Mirrlees99}
Consider the Mirrlees' problem. Note that the solution of Mirrlees' problem does not change if we add the constraint $y\in [-2,2]$ into the problem.
 \begin{eqnarray*}
~~~~~\min && (x-2)^2+(y-1)^2\\
{\rm s.t.} && y \in S(x),
\end{eqnarray*}
where $S(x)$ is the solution set of the lower level program
 \begin{eqnarray*}
~~~~~\min && -x \exp[-(y+1)^2]-\exp[-(y-1)^2]\\
{\rm s.t.} && y\in [-2,2].
\end{eqnarray*}
It was shown in \cite{Mirrlees99} that the unique optimal solution is $(\bar{x},\bar{y})$ with $\bar x=1,$ $\bar y\approx 0.958$  being the positive solution of the equation
$$(1+y) =(1-y) \exp[4y].$$
\end{eg}
In our test, we chose the initial point $(x_0,y_0)=(0.5,0.3)$ and the parameters $\beta=0.8,\ \sigma_1=\sigma_2=10^{-6},\  \rho_0=100,\ r_0=100,\ \hat{\eta}=5*10^5,\ \sigma=10,\sigma'=10$ and $\varepsilon=7*10^{-5}$, $\varepsilon'=10^{-8}$, $\epsilon_1=10^{-6}$.  Since the stopping criterion $\|(x^{k-1},y^{k-1})-(x^k,y^k)\|\leq \epsilon_1$ hold, we terminate at the 8th iteration with   $(x^k,y^k)=(1,0.95759)$. It seems that the sequence converges to  $(\bar x,\bar y)$.

 Since
\begin{eqnarray*}
&&\nabla f(x^k,y^k)-(\nabla \gamma_{\rho_k}(x^{k}),0)=(0.97665 ,0.00015) ,\\
&& \nabla (\nabla_y f)(x^k,y^k)=(0.084814, 1.70047) ,
\end{eqnarray*}
by virtue of the continuity of the the gradients it is easy to see that the vectors
$$\nabla f(\bar x,\bar y)-(\displaystyle\lim_{k\rightarrow \infty} \nabla \gamma_{\rho_k}(x^{k}),0) \mbox{ and } \nabla (\nabla_y f)(\bar x,\bar y)$$
  are linearly independent. 
 Thus the WNNAMCQ holds at $(\bar x,\bar y)$ and our algorithm guarantees that $(\bar x,\bar y)$ is a stationary point of (CP). Indeed, $(\bar x,\bar y)$ is  the unique global minimizer of the Mirrlees' problem.

\begin{eg}\label{test3.14} \rm\cite[Example 3.14]{test} The bilevel program
\begin{eqnarray*}
\min && F(x,y):=(x-\frac{1}{4})^2+y^2\\
{\rm s.t.}
&& y\in S(x):=\argmin\limits_{y\in [-1,1]} f(x,y):={\textstyle \frac{y^3}{3}-x y }
\end{eqnarray*}
\end{eg}
has the optimal solution point $(\bar{x},\bar{y})=(\frac{1}{4},\frac{1}{2})$ with an objective value of $\frac{1}{4}$.

In our test, we chose the initial point $(x_0,y_0)=(0.3,0.3)$ and the parameters $\beta=0.9,\ \sigma_1=\sigma_2=10^{-6},\  \rho_0=100,\ r_0=100,\ \hat{\eta}=5000,\ \sigma=10,\sigma'=10$ and $\varepsilon=5*10^{-6}$, $\varepsilon'=10^{-8}$, $\epsilon_1=5*10^{-6}$.  Since the stopping criterion $\|(x^{k-1},y^{k-1})-(x^k,y^k)\|\leq \epsilon_1$ hold, we terminate at the 7th iteration and obtain an point  $(x^k,y^k)=(0.25,0.5)$.  It seems that the sequence converges to  $(\bar x,\bar y)$.

  Since
\begin{eqnarray*}
&&\nabla f(x^k,y^k)-(\nabla \gamma_{\rho_k}(x^{k}),0)=(-1.5,2.3*10^{-7}),\\
&& \nabla (\nabla_y f)(x^k,y^k)=(-1,1),
\end{eqnarray*}
by virtue of the continuity of the the gradients it is easy to see that the vectors
$$\nabla f(\bar x,\bar y)-(\displaystyle\lim_{k\rightarrow \infty} \nabla \gamma_{\rho_k}(x^{k}),0) \mbox{ and } \nabla (\nabla_y f)(\bar x,\bar y)$$
  are linearly independent.  
 Thus the WNNAMCQ holds at $(\bar x,\bar y)$ and our algorithm guarantees that $(\bar x,\bar y)$ is a stationary point of (CP). Indeed, $(\bar x,\bar y)$ is  the unique global minimizer of the  problem.

\begin{eg}\label{test3.20} \rm\cite[Example 3.20]{test} The bilevel program
\begin{eqnarray*}
\min && F(x,y):=(x-0.25)^2+ y^2\\
{\rm s.t.}
&& y\in S(x):=\argmin\limits_{y\in [-1,1]} f(x,y):={\textstyle  \frac{1}{3}y^3-x^2 y }
\end{eqnarray*}
\end{eg}
has the optimal solution point $(\bar{x},\bar{y})=(\frac{1}{2},\frac{1}{2})$ with an objective value of $\frac{5}{16}$.

In our test, we chose the parameters $\beta=0.9,\ \sigma_1=\sigma_2=10^{-6},\  \rho_0=100,\ r_0=100,\ \hat{\eta}=500,\ \sigma=10,\sigma'=10$ and $\varepsilon=10^{-6}$, $\varepsilon'=10^{-8}$, $\epsilon_1=10^{-6}$.
We chose the initial point $(x_0,y_0)=(0.3,0.3)$.  Since the stopping criterion $\|(x^{k-1},y^{k-1})-(x^k,y^k)\|\leq \epsilon_1$ hold, we terminate at the 8th iteration and obtain an point  $(x^k,y^k)=(0.4999998,0.4999998)$.  It seems that the sequence converges to  $(\bar x,\bar {y})$.

  Since
\begin{eqnarray*}
&&\nabla f(x^k,y^k)-(\nabla \gamma_{\rho_k}(x^{k}),0)=(-1.49989,0),\\
&& \nabla (\nabla_y f)(x^k,y^k)=(-1,1),
\end{eqnarray*}
 by virtue of the continuity of the the gradients it is easy to see that the vectors
$$\nabla f(\bar x,\bar {y}_1)-(\displaystyle\lim_{k\rightarrow \infty} \nabla \gamma_{\rho_k}(x^{k}),0) \mbox{ and } \nabla (\nabla_y f)(\bar x,\bar{y}_1)$$
  are linearly independent.  
 Thus the WNNAMCQ holds at $(\bar x,\bar {y})$ and our algorithm guarantees that $(\bar x,\bar y)$ is a stationary point of (CP). Indeed, $(\bar x,\bar y)$ is  the unique global minimizer of the  problem.

\section{Conclusions}
In this paper, we propose a smoothing SQP method for solving nonsmooth and nonconvex  optimization problems with Lipschitz inequality and equality constraints. The algorithm is applicable even to  degenerate  constrained optimization problems  which do not satisfy the  GMFCQ, the standard  constraint qualification for a local minimizer to satisfy the KKT conditions. Our main motivation comes from solving the bilevel program which is nonsmooth, nonconvex and never satisfies 	 the GMFCQ.  In this paper, we have proposed the concept of the WGMFCQ (equivalently WNNAMCQ),  a weaker version of the GMFCQ,  and have shown  
 the global convergence of the smoothing SQP algorithm under the WGMFCQ. Moreover we have demonstrated the applicability of the smoothing SQP algorithm for solving the combined program of a simple bilevel program with a nonconvex  lower level problem. For smooth optimization problem, it is well-known that the SQP methods converge very fast  when the iterates are close to the solution. The rapid local convergence of the SQP is due to the fact that the positive definite matrix $W_k$ in the SQP subproblem  is an approximation of the Hessian matrix of the Lagrangian function.   For our nonsmooth problem, the Lagrangian function is only locally Lipschitz and no classical Hessian matrix can be defined. However it would be interesting to study the local behaviour of the smoothing SQP algorithm 
 by using the generalized second order subderivatives (\cite{var}) of the Lagrangian function. This  remains a topic of our future research.

\end{document}